\documentclass[11pt,reqno]{amsart}
\usepackage[T1]{fontenc}
\usepackage[utf8]{inputenc}
\usepackage{lmodern}
\usepackage[a4paper,margin=28mm]{geometry}
\usepackage{amssymb}
\usepackage[expansion=false]{microtype}
\usepackage{xcolor}
\usepackage[colorlinks=true,linkcolor=blue!45!black,citecolor=blue!45!black,urlcolor=blue!45!black]{hyperref}
\hypersetup{pdftitle={Wirsing type approximation over totally real and totally p-adic extensions}}
\numberwithin{equation}{section}
\newtheorem{theorem}{Theorem}[section]
\newtheorem{proposition}[theorem]{Proposition}
\newtheorem{lemma}[theorem]{Lemma}

\theoremstyle{definition}

\theoremstyle{remark}
\newtheorem{remark}[theorem]{Remark}
\newcommand{\Q}{\mathbb Q}
\newcommand{\Z}{\mathbb Z}
\newcommand{\R}{\mathbb R}
\newcommand{\Qp}{\mathbb Q_p}
\newcommand{\Zp}{\mathbb Z_p}
\newcommand{\hl}{\widehat\lambda_{n,p}(\xi)}
\title[Wirsing type approximation over extensions]{Wirsing type approximation over totally real and totally $p$-adic extensions}
\author{Anup B. Dixit}
\address{Department of Mathematics, The Institute of Mathematical Sciences (HBNI), CIT Campus, IV Cross Road, Chennai, India--600113}
\email{anupdixit@imsc.res.in}
\author{Hitesh Kumar}
\address{Department of Mathematics, The Institute of Mathematical Sciences (HBNI), CIT Campus, IV Cross Road, Chennai, India--600113}
\email{hiteshk@imsc.res.in}
\date{}
\begin{document}

\begin{abstract}
In this paper, we establish Wirsing-type bounds for approximation of real and $p$-adic transcendental numbers by totally real and totally $p$-adic algebraic numbers, respectively. More precisely, for every $n\ge 1$, we show that every transcendental $\xi\in\Q_p$ admits infinitely many totally $p$-adic algebraic approximants $\alpha$ of degree at most $n$ satisfying $0<|\xi-\alpha|_p\ll H(\alpha)^{-1-\lceil n/2\rceil}$. We improve this bound for even $n\ge 4$. For each $n\geq 2$, we also show that every real number $\xi$ that is not algebraic of degree at most $n$ admits totally real approximants of degree at most $n$ with a bound stronger than Wirsing's. In both cases, we can also restrict the approximants to be totally $\ell$-adic for every prime $\ell$ in a fixed finite set.
\end{abstract}
\maketitle

\section{Introduction}

For a nonzero polynomial $P(X)=a_0+a_1X+\cdots+a_dX^d\in\Z[X]$ of degree $d$, we define its height by
$$
H(P)=\max_{0\le i\le d}|a_i|
$$
and set $H(0)=0$. For an algebraic number $\alpha$, let $H(\alpha)$ be the height of its primitive minimal polynomial over $\Z$. A classical theorem of Wirsing \cite{Wirsing} states that, for every integer $n\ge1$ and every real number $\xi$ which is not algebraic of degree at most $n$, there are infinitely many real algebraic numbers $\alpha$ of degree at most $n$ such that
$$
0<|\xi-\alpha|\ll H(\alpha)^{-(n+3)/2}.
$$

Wirsing conjectured that the exponent $(n+3)/2$ can be replaced by $n+1-\varepsilon$, for every $\varepsilon>0$. The case $n=2$ was proved by Davenport--Schmidt \cite{DSquadratic}. Subsequent improvements in general degree were obtained by Bernik--Tishchenko \cite{BT} and Tishchenko \cite{Tishchenko,Tsishchanka}. We refer to \cite[Chapter~3]{Bugeaud} for a comprehensive account of the classical developments in Wirsing's problem.\\

More recently, Badziahin--Schleischitz \cite{BS} obtained the error exponent $1+n/\sqrt3-\varepsilon$ for $n\ge4$. Po\"els \cite{Poels} subsequently improved this to $1+n/(2-\log 2)-\varepsilon$ for $n\ge2$. These bounds hold for every transcendental real number $\xi$ and every $\varepsilon>0$.\\

The corresponding problem over $\Qp$ was studied by Morrison \cite{Morrison} and Teuli\'e \cite{Teulie} (also see \cite[Chapter~9]{Bugeaud}). Here both $\xi$ and the algebraic approximants $\alpha$ lie in $\Qp$ and the error is measured by $|\xi-\alpha|_p$. In \cite{Dixit}, the first author showed that one may obtain the error exponent $n/(2-\log 2)-\varepsilon$, providing a $p$-adic counterpart of Po\"els's bound.\\

In this paper, we initiate the study of Wirsing's problem in which the approximants are restricted to an \emph{infinite algebraic extension} $K\subset\overline{\Q}$, instead of allowing approximants from the whole of $\overline{\Q}$. Over such an extension $K$, the quality of approximation need not improve indefinitely as the degree increases. For example, let $\Q^{(2)}$ be the compositum of all quadratic extensions of $\Q$. For $n\ge2$ and a transcendental real number $\xi$, let $\omega_n(\xi)$ be the supremum of the real numbers $w$ for which
\begin{equation*}
0<|\xi-\alpha|<H(\alpha)^{-w-1}
\end{equation*}
holds for infinitely many distinct $\alpha\in\Q^{(2)}\cap\R$ of degree at most $n$. As shown in Proposition~\ref{prop:multiquadratic-obstruction}, for almost every real $\xi$ 
\begin{equation*}
\omega_n(\xi)\le5-\frac{4}{n}.
\end{equation*}
Thus, for almost all $\xi\in\R$, $\omega_n(\xi) < 5$ and hence Wirsing-type approximation fails over $\Q^{(2)}$.\\

We consider approximations by totally real numbers $\Q^{tr}$, totally $p$-adic numbers $\Q^{tp}$ and finite intersections of totally $\ell$-adic fields as $\ell$ ranges over a finite set of primes. For each of these infinite extensions, we establish approximations for transcendental numbers of a similar quality as Wirsing's.\\

Let $\ell$ be a prime. Recall that an algebraic number $\alpha$ is \emph{totally $\ell$-adic} if its minimal polynomial splits completely over $\Q_\ell$. Equivalently, if $L=\Q(\alpha)$ and $d=[L:\Q]$, then
$$
L\otimes_{\Q}\Q_\ell\simeq\Q_\ell^d.
$$
Thus total $\ell$-adicity requires every conjugate to lie in $\Q_\ell$, whereas the condition $\alpha\in\Q_\ell$ requires only one such embedding. Similarly, $\alpha$ is \emph{totally real} if all its conjugates are real.\\

Throughout the paper, $S$ is a fixed finite set of primes, possibly empty. We say that $\alpha$ is totally $S$-adic if it is totally $\ell$-adic for every $\ell\in S$. Our first result concerns $p$-adic approximation.

\begin{theorem}\label{thm:basic}
Let $p$ be a prime, $n\ge1$, and $\xi\in\Qp$ transcendental. There is a constant $C=C(n,p,\xi,S)>0$ such that infinitely many distinct algebraic numbers $\alpha\in\Qp$ of degree at most $n$ are totally $(S\cup\{p\})$-adic and satisfy
\begin{equation}\label{eq:basic}
0<|\xi-\alpha|_p\le C H(\alpha)^{-1-\lceil n/2\rceil}.
\end{equation}
\end{theorem}

Taking $S=\varnothing$ gives $p$-adic approximation by totally $p$-adic numbers. Taking $S=\{\ell\}$, with $\ell\ne p$, gives $p$-adic approximation by numbers which are simultaneously totally $p$-adic and totally $\ell$-adic.\\

For real approximation, we obtain a stronger bound. Let
\begin{equation}\label{eq:Bn}
B_n :=\frac n2+\frac{1-\log 2}{2}\sqrt n+\frac13\quad(n\ge2).
\end{equation}

\begin{theorem}\label{thm:real}
Let $n\ge2$ and let $\xi\in\R$ be not algebraic of degree at most $n$. For every $\varepsilon>0$, infinitely many distinct algebraic numbers $\alpha\in\R$ of degree at most $n$ are both totally real and totally $S$-adic, and satisfy
\begin{equation*}
0<|\xi-\alpha|\le H(\alpha)^{-1-B_n+\varepsilon}.
\end{equation*}
\end{theorem}

Consequently, for every $n\ge1$ and every real $\xi$ which is not algebraic of degree at most $n$, infinitely many totally real, totally $S$-adic numbers $\alpha$ of degree at most $n$ satisfy
\begin{equation*}
0<|\xi-\alpha|\le C H(\alpha)^{-(n+3)/2}
\end{equation*}
for a constant $C=C(n,\xi,S)>0$.
This is because for $n\ge2$
\begin{equation*}
B_n-\frac{n+1}{2}=\frac{1-\log 2}{2}\sqrt n-\frac16>0.
\end{equation*}

\subsection{Method of proof} Our proofs use the inhomogeneous transference method originating in the work of Davenport--Schmidt \cite{DS}; see also \cite[Proposition~3.3]{Bugeaud}. The $p$-adic version of this method was developed by Teuli\'e \cite{Teulie}. We combine transference with local conditions on the approximating polynomials to ensure the prescribed splitting. This reduces the approximation problem to bounds for uniform simultaneous approximation, from which our main estimates follow. In even degrees, we sharpen the $p$-adic bound by using a method of Schleischitz \cite{Schleischitz}.

\bigskip

\section{Preliminaries}\label{sec:prelim}

In this section, we develop the ingredients needed to prove our main theorems.\\

Write $A\asymp B$ when both $A\ll B$ and $B\ll A$ hold. Recall the following standard estimates on heights (see \cite[Lemmas~A.3 and~A.4]{Bugeaud}).

\begin{lemma}\label{lem:height}
If $P\in\Z[X]$ is nonzero and has degree at most $n$, and $P(\alpha)=0$, then
\begin{equation*}
H(\alpha)\le 2^n\sqrt{n+1}\,H(P).
\end{equation*}
For fixed $c\in\Q^{\times}$ and $\deg\alpha\le n$, one also has
\begin{equation*}
H(c\alpha)\asymp_{c,n}H(\alpha).
\end{equation*}
Multiplication by $c$ preserves degree, total reality, and total $S$-adicity for every finite set $S$.
\end{lemma}
\medskip

\subsection{Lattice theory}\label{subsec:lattices}

A full-rank lattice in $\R^d$ is a group $L=B\Z^d$ with $B$ an invertible real matrix. We denote by $B^t$ the transpose of $B$. Its determinant is $\det L:=|\det B|$ and its dual is
\begin{equation*}
L^*:=\{y\in\R^d:y\cdot x\in\Z\text{ for every }x\in L\}
=B^{-t}\Z^d.
\end{equation*}
Thus $\det L^*=(\det L)^{-1}$. Write $\|x\|_2$ and $\|x\|_\infty$ for the Euclidean and maximum norms. They satisfy
\begin{equation*}
\|x\|_\infty\le\|x\|_2\le\sqrt d\,\|x\|_\infty.
\end{equation*}
Denote the first Euclidean minimum by
\begin{equation*}
\mu_1(L):=\min\{\|x\|_2:x\in L\setminus\{0\}\}
\end{equation*}
and the covering radius of $L$ by
\begin{equation*}
\rho(L):=\sup_{u\in\R^d}\inf_{a\in L}\|u-a\|_2.
\end{equation*}
We shall use the following transference bound of Banaszczyk {\cite[Theorem~2.2]{Banaszczyk}}.

\begin{lemma}[Banaszczyk]\label{lem:cover}
For every full-rank lattice $L\subset\R^d$, we have
\begin{equation*}
\rho(L)\le \frac{d}{2\mu_1(L^*)}.
\end{equation*}
\end{lemma}
\medskip

\subsection{Uniform exponent}\label{subsec:uniform}

Let $v$ denote a place over $\Q$, either $\infty$ or a rational prime $p$. Let $n\ge1$ and assume that $\xi$ is not algebraic of degree at most $n$. Equivalently, $1,\xi,\ldots,\xi^n$ are linearly independent over $\Q$. For $\xi\in\R$, define $\widehat\lambda_{n,\infty}(\xi)$ as the supremum of the real numbers $\lambda$ such that, for every sufficiently large real $T$, there is a nonzero vector $(x_0,\ldots,x_n)\in\Z^{n+1}$ with
\begin{equation}\label{eq:realuniform}
|x_0|\le T,\qquad
\max_{1\le i\le n}|x_0\xi^i-x_i|\le T^{-\lambda}.
\end{equation}
We follow the convention of \cite{PRreal}. Bounding all the coordinates by $T$ gives the same exponent, since the above inequalities imply $|x_i|\ll_\xi T$.\\

For $\xi\in\Qp$ satisfying the same independence condition, define $\hl$ as the supremum of all real numbers $\lambda$ such that, for every sufficiently large real $T$, there is a nonzero vector $(x_0,\ldots,x_n)\in\Z^{n+1}$ satisfying
\begin{equation}\label{eq:uniform}
\max_{0\le i\le n}|x_i|\le T,\qquad
0<\max_{1\le i\le n}|x_0\xi^i-x_i|_p\le T^{-1-\lambda}.
\end{equation}
The vectors in this definition are not required to be primitive. This is the normalization in \cite[Definition~2.8]{BBS}.\\

For $\lambda>\widehat\lambda_{n,\infty}(\xi)$, the system \eqref{eq:realuniform} (likewise for $\lambda>\widehat\lambda_{n,p}(\xi)$, the system \eqref{eq:uniform}) has no nonzero integral solution for arbitrarily large $T$.  

\begin{lemma}\label{lem:scaling}
For $v=\infty$ or $v=p$, a number $\xi\in\Q_v$ which is not algebraic of degree at most $n$ and $c\in\Q^\times$, we have
\begin{equation*}
\widehat\lambda_{n,v}(c\xi)=\widehat\lambda_{n,v}(\xi).
\end{equation*}
\end{lemma}
\begin{proof}
For $c=a/b$, with nonzero integers $a,b$, the map
\begin{equation*}
y_i:=a^i b^{n-i}x_i\quad(0\le i\le n)
\end{equation*}
preserves integrality and gives
\begin{equation*}
y_0(c\xi)^i-y_i=a^i b^{n-i}(x_0\xi^i-x_i).
\end{equation*}
Thus heights and errors change by fixed factors. Applying the same argument to $c^{-1}$ proves the invariance.
\end{proof}

We also recall the following two bounds for uniform exponents.

\begin{theorem}[Badziahin--Bugeaud--Schleischitz]\label{thm:BBS}
Let $n\ge1$, let $\xi\in\Qp$ be transcendental, and put $k:=\lceil n/2\rceil$. Then
\begin{equation*}
\frac1n\le\hl\le\frac1k.
\end{equation*}
\end{theorem}

The lower bound follows from Dirichlet's theorem (\cite[equation~(2.4)]{BBS}). The upper bound follows from \cite[Theorem~2.10]{BBS}. Therefore, there is a constant $c=c(n,p,\xi)>0$ such that, for arbitrarily large real $T$, the system
\begin{equation*}
0<\max_{0\le i\le n}|x_i|\le T,
\qquad
\max_{1\le i\le n}|x_0\xi^i-x_i|_p\le cT^{-1-1/k}
\end{equation*}
has no integral solution. We shall use this to obtain Theorem~\ref{thm:basic} without a loss in the exponent.\\

For real $\xi$, we use the following result of Po\"els--Roy {\cite[Theorem~1.1]{PRreal}}.

\begin{theorem}[Po\"els--Roy]\label{thm:PRreal}
For every $n\ge2$ and every real $\xi$ which is not algebraic of degree at most $n$,
\begin{equation*}
\frac1n\le\widehat\lambda_{n,\infty}(\xi)\le\frac1{B_n},
\end{equation*}
where $B_n$ is defined in \eqref{eq:Bn}.
\end{theorem}

Again the lower bound follows from Dirichlet's theorem. The upper bound gives arbitrarily large parameters for which the system defining $\widehat\lambda_{n,\infty}(\xi)$ has no nonzero solution whenever $\lambda>B_n^{-1}$. This accounts for the $\varepsilon$ in Theorem~\ref{thm:real}.
\medskip

\subsection{Congruences and complete splitting}\label{subsec:local}

For a prime $\ell$, normalize $|\ell|_\ell=\ell^{-1}$ and write $v_\ell(\ell)=1$, with $v_\ell(0):=+\infty$. We recall the strong Hensel's lemma (see \cite[Theorem~4.1]{Conrad}).

\begin{lemma}[Hensel's lemma]\label{lem:hensel}
Let $G\in\Z_\ell[X]$ and $a\in\Z_\ell$. If
\begin{equation*}
v_\ell(G(a))>2v_\ell(G'(a)),
\end{equation*}
then $G$ has a unique root $\beta$ in the ball
\begin{equation*}
v_\ell(\beta-a)>v_\ell(G'(a)).
\end{equation*}
If $G(a)\ne0$, then
\begin{equation*}
v_\ell(\beta-a)=v_\ell(G(a))-v_\ell(G'(a)).
\end{equation*}
\end{lemma}
\medskip

For $n\ge1$ and a prime $\ell$, put
\begin{equation*}
E_\ell:=v_\ell((n-1)!).
\end{equation*}
Given an integer $r$, define
\begin{equation}\label{eq:Gell}
G_r(X):=\prod_{i=0}^{n-1}(X-r-i).
\end{equation}

\begin{proposition}\label{prop:split}
If $P\in\Z[X]$ has degree at most $n$ and
\begin{equation*}
P\equiv G_r\pmod{\ell^{2E_\ell+1}}
\end{equation*}
coefficientwise, then $P$ has degree $n$ and splits into $n$ distinct linear factors over $\Q_\ell$. Its roots can be labelled $\beta_0,\ldots,\beta_{n-1}$ so that
\begin{equation*}
v_\ell(\beta_i-r-i)\ge E_\ell+1.
\end{equation*}
If in addition $\ell=p$ and $\xi-r\in p^{E_p+1}\Zp$, then
\begin{equation*}
|P'(\xi)|_p=|(n-1)!|_p,
\qquad
|\xi-\beta_0|_p=\frac{|P(\xi)|_p}{|(n-1)!|_p}.
\end{equation*}
\end{proposition}
\begin{proof}
The leading coefficient of $P$ is $1$ modulo $\ell^{2E_\ell+1}$, so it is a unit in $\Z_\ell$ and $\deg P=n$. For $0\le i\le n-1$, set
\begin{equation*}
e_i:=v_\ell\bigl(i!(n-1-i)!\bigr).
\end{equation*}
Since $\binom{n-1}{i}$ is an integer, $e_i\le E_\ell$. At $r+i$ we have
\begin{equation*}
P(r+i)\in\ell^{2E_\ell+1}\Z_\ell,\qquad
P'(r+i)\equiv(-1)^{n-1-i}i!(n-1-i)!
\pmod{\ell^{2E_\ell+1}\Z_\ell}.
\end{equation*}
The derivative has valuation $e_i$ and $2E_\ell+1>2e_i$. By Hensel's lemma (Lemma~\ref{lem:hensel}), we have a root $\beta_i$ with
\begin{equation*}
v_\ell(\beta_i-r-i)\ge 2E_\ell+1-e_i\ge E_\ell+1.
\end{equation*}
If $i\ne j$, then $1\le|i-j|\le n-1$, so $v_\ell(i-j)\le E_\ell$ and the ultrametric inequality gives
\begin{equation*}
v_\ell(\beta_i-\beta_j)=v_\ell(i-j).
\end{equation*}
Thus these are $n$ distinct roots of $P$, and since $P$ has degree $n$, they are all its roots. Note that differentiation preserves coefficientwise congruence modulo $p^{2E_p+1}$. Since $\xi-r\in p^{E_p+1}\Zp$ and $2E_p+1\ge E_p+1$,
\begin{equation*}
P'(\xi)\equiv G_r'(r)=(-1)^{n-1}(n-1)!
\pmod{p^{E_p+1}\Zp}.
\end{equation*}
This proves the derivative identity. For $1\le i\le n-1$, the root estimate in Proposition~\ref{prop:split} implies
$v_p(\xi-\beta_i)=v_p(i)$.
Now, factoring $P$ over $\Qp$, its leading coefficient is a unit, and consequently
\begin{equation*}
|P(\xi)|_p=|\xi-\beta_0|_p\prod_{i=1}^{n-1}|i|_p
=|\xi-\beta_0|_p\,|(n-1)!|_p.
\end{equation*}
\end{proof}

\begin{remark}
The proposition also applies when $n>\ell$, since its root neighborhoods are disjoint even if some centers agree modulo $\ell$. 
\end{remark}

\begin{lemma}\label{lem:CRT}
Let $S$ be a finite set of primes and choose an integer $r_\ell$ for each $\ell\in S$. Put
\begin{equation*}
D:=\prod_{\ell\in S}\ell^{2E_\ell+1}.
\end{equation*}
There is a fixed monic $F\in\Z[X]$ of degree $n$ with
\begin{equation*}
F\equiv G_{r_\ell}\pmod{\ell^{2E_\ell+1}}\qquad(\ell\in S),
\end{equation*}
where $G_{r_\ell}$ is as in \eqref{eq:Gell}, such that every $P\in\Z[X]$ of degree at most $n$ satisfying $P\equiv F\pmod D$ splits completely over each $\Q_\ell$, $\ell\in S$. If $S$ is nonempty, such a $P$ has degree exactly $n$.
\end{lemma}
\begin{proof}
By definition, all $G_{r_\ell}$ are monic of degree $n$. For each $0\le i<n$, the Chinese remainder theorem allows us to choose an integer coefficient of $X^i$ in $F$ congruent to the coefficient of $X^i$ in $G_{r_\ell}$ modulo $\ell^{2E_\ell+1}$ for every $\ell\in S$. Taking the leading coefficient of $F$ to be $1$ gives the required monic polynomial. Then, applying Proposition~\ref{prop:split} at each prime in $S$, we obtain the desired conclusion.
\end{proof}
\medskip

\subsection{Different approximation behaviour in the multiquadratic extension}\label{subsec:multiquadratic}

Recall that, for $n\ge2$ and a transcendental real number $\xi$, the exponent $\omega_n(\xi)$ is the supremum of the real numbers $w$ for which
\begin{equation*}
0<|\xi-\alpha|<H(\alpha)^{-w-1}
\end{equation*}
holds for infinitely many distinct $\alpha\in\Q^{(2)}\cap\R$ of degree at most $n$.

\begin{proposition}\label{prop:multiquadratic-obstruction}
For almost every real $\xi$ and any integer $n\ge2$,
\begin{equation*}
\omega_n(\xi)\le5-\frac4n.
\end{equation*}
\end{proposition}
\begin{proof}
We first count the possible fields generated by the approximants. Every finite subfield of $\Q^{(2)}$ is abelian. For an abelian field $L/\Q$ of degree $m\ge2$, let $f_L$ be its conductor, $D_L$ its absolute discriminant and $\mathcal X_L$ its group of Dirichlet characters. By the conductor--discriminant formula, we have
\begin{equation*}
D_L=\prod_{\chi\in\mathcal X_L}f_\chi
\,\,\,\text{and}\,\,\, f_L=\mathop{\mathrm{lcm}}_{\chi\in\mathcal X_L}f_\chi.
\end{equation*}
 If $\ell^a\Vert f_L$, the characters with $v_\ell(f_\chi)<a$ form a proper subgroup of $\mathcal X_L$. At least half of its $m$ characters therefore have $v_\ell(f_\chi)=a$. Applying this at every prime dividing $f_L$ yields
\begin{equation*}
D_L\ge f_L^{m/2}.
\end{equation*}

For fixed $m$, the number of characters modulo $f$ of order dividing $m$ is at most $m^{\omega(f)+1}$, where $\omega(f)$ counts the distinct prime factors of $f$. Indeed, the unit group modulo an odd prime power is cyclic, while at a power of $2$ it is a product of at most two cyclic groups. A group of $m$ characters is specified by a list of $m$ such characters, so there are at most
\begin{equation*}
m^{m(\omega(f)+1)}\ll_{m,\eta} f^\eta
\end{equation*}
possible degree-$m$ abelian fields of conductor $f$, for every $\eta>0$. Since $D_L\le X$ implies $f_L\le X^{2/m}$, summing over the conductors gives
\begin{equation}\label{eq:abelian-field-count}
\#\{L/\Q\text{ abelian}:[L:\Q]=m,\ D_L\le X\}
\ll_{m,\varepsilon}X^{2/m+\varepsilon}.
\end{equation}

We next count generators in each field. For a degree-$m$ algebraic number $\alpha$, define its Mahler measure by
\begin{equation*}
M(\alpha):=|a_m|\prod_{j=1}^{m}\max\{1,|\alpha_j|\},
\end{equation*}
where $a_m$ is the leading coefficient of its primitive minimal polynomial and $\alpha_1,\ldots,\alpha_m$ are its roots. The standard height comparison gives
\begin{equation*}
M(\alpha)\le\sqrt{m+1}\,H(\alpha).
\end{equation*}
If distinct real numbers $\alpha,\beta$ both generate the same degree-$m$ field and have height at most $B$, the product formula gives
\begin{equation*}
1\le2^m|\alpha-\beta|M(\alpha)M(\beta)
\le2^m(m+1)|\alpha-\beta|B^2.
\end{equation*}
Thus they are separated by at least $2^{-m}(m+1)^{-1}B^{-2}$. In a fixed bounded interval $I$, each field has only $O_{m,I}(B^2)$ such generators, with a constant independent of the field.

For a generator $\alpha$ of height at most $B$, Silverman's discriminant--height inequality \cite{Silverman} gives
\begin{equation*}
D_{\Q(\alpha)}\ll_m H(\alpha)^{2m-2}\le B^{2m-2}.
\end{equation*}
Combining this with \eqref{eq:abelian-field-count} gives
\begin{equation}\label{eq:multiquadratic-count}
\#\{\alpha\in\Q^{(2)}\cap I:\deg\alpha=m,\ H(\alpha)\le B\}
\ll_{m,I,\varepsilon}B^{6-4/m+\varepsilon}.
\end{equation}
Rationals contribute $O_I(B^2)$. Since $6-4/m$ increases with $m$, summing over degrees at most $n$ gives the same bound with exponent $6-4/n+\varepsilon$. Fix $n\ge2$ and $\tau>6-4/n$, and choose $0<\varepsilon<\tau-(6-4/n)$. For approximants with $2^j\le H(\alpha)<2^{j+1}$, the intervals
\begin{equation*}
|\xi-\alpha|<H(\alpha)^{-\tau}
\end{equation*}
have total length
\begin{equation*}
\ll_{n,I,\varepsilon}2^{j(6-4/n+\varepsilon-\tau)}.
\end{equation*}
We now apply the Borel--Cantelli lemma. Fix an integer $R\ge1$ and restrict the targets to $[-R,R]$. Since $H(\alpha)^{-\tau}\le1$, any approximant to such a target lies in $I=[-R-1,R+1]$. For each integer $j\ge0$, let $A_j$ be the union, over all $\alpha\in\Q^{(2)}\cap I$ with $\deg\alpha\le n$ and $2^j\le H(\alpha)<2^{j+1}$, of the intervals
\begin{equation*}
(\alpha-H(\alpha)^{-\tau},\alpha+H(\alpha)^{-\tau})\cap[-R,R].
\end{equation*}
Writing $\operatorname{meas}$ for Lebesgue measure, the preceding estimate gives
\begin{equation*}
\sum_{j=0}^{\infty}\operatorname{meas}(A_j)
\ll_{n,R,\varepsilon}\sum_{j=0}^{\infty}
2^{j(6-4/n+\varepsilon-\tau)}<\infty.
\end{equation*}
The Borel--Cantelli lemma therefore implies that $\operatorname{meas}(\limsup_{j\to\infty}A_j)=0$. Thus almost every $\xi\in[-R,R]$ belongs to only finitely many $A_j$. Each height range contains only finitely many algebraic numbers of degree at most $n$, so such a $\xi$ admits only finitely many approximants satisfying $|\xi-\alpha|<H(\alpha)^{-\tau}$. 

Taking the union of these null exceptional sets over integers $R\ge1$, integers $n\ge2$ and rational $\tau>6-4/n$ and excluding the algebraic real numbers, leaves a set of full measure in $\R$. On this set, $\omega_n(\xi)\le5-4/n$ for every $n\ge2$. Otherwise, infinitely many approximants would satisfy the inequality with exponent $w+1>6-4/n$, contradicting the finiteness assertion for any rational $\tau$ between $6-4/n$ and $w+1$.
\end{proof}

\section{Algebraic approximation with prescribed splitting}\label{sec:transference}

The key step in proving our main results is the following transference theorem, which relates approximation to $\xi$ by algebraic numbers satisfying prescribed splitting conditions to uniform simultaneous rational approximation to its first $n$ powers.

\begin{theorem}\label{thm:reciprocal}
Let $v$ be either $\infty$ or a prime $p$, let $n\ge1$ and let $\xi\in\Q_v$ be not algebraic of degree at most $n$. If $\lambda>\widehat\lambda_{n,v}(\xi)$, then infinitely many distinct algebraic numbers $\alpha\in\Q_v$, of degree at most $n$ and totally $S$-adic, satisfy
\begin{equation}\label{eq:reciprocal}
0<|\xi-\alpha|_v\ll H(\alpha)^{-1-1/\lambda}.
\end{equation}
When $v=p$, the approximants can also be required to be totally $p$-adic. When $v=\infty$, they can also be required to be totally real. The implied constant may depend on $n,v,\xi,S,$ and $\lambda$.
\end{theorem}

Together with Theorem~\ref{thm:BBS} and its nonexistence assertion, this gives Theorem~\ref{thm:basic}, including the fixed multiplicative constant. Combining it with Theorem~\ref{thm:PRreal} gives Theorem~\ref{thm:real}. We prove the $p$-adic and real cases in Subsections~\ref{subsec:padicproof} and~\ref{subsec:realproof}, respectively.

\subsection{The \texorpdfstring{$p$}{p}-adic case}\label{subsec:padicproof}

For the $p$-adic proofs, it suffices to treat $\xi\in\Zp$. Indeed, choose an integer $r\ge0$ with $\zeta=p^r\xi\in\Zp$. By Lemma~\ref{lem:scaling}, we have $\widehat\lambda_{n,p}(\zeta)=\widehat\lambda_{n,p}(\xi)$. Given approximants $\beta$ to $\zeta$, set $\alpha=p^{-r}\beta$. Lemma~\ref{lem:height} preserves their degree and splitting conditions and gives $H(\alpha)\asymp H(\beta)$, while
\begin{equation*}
|\xi-\alpha|_p=p^r|\zeta-\beta|_p.
\end{equation*}
Thus the approximation exponent is unchanged. The condition that $\xi$ is not algebraic of degree at most $n$ is also preserved by this scaling.\\

Assume $\xi\in\Zp$. Fix $q=p^h$, and choose integers $z_i$ representing $\xi^i$ modulo $q\Zp$, with $z_0=1$. Put $z:=(z_0,\ldots,z_n)$ and
\begin{equation*}
L:=\{a\in\Z^{n+1}:a\cdot z\equiv0\pmod q\}.
\end{equation*}
Under the identification of polynomials and coefficient vectors, this is the lattice of $Q\in\Z[X]$ of degree at most $n$ for which $Q(\xi)\in q\Zp$.

\begin{lemma}\label{lem:dual}
The lattice $L$ has determinant $q$ and dual
\begin{equation*}
L^*=\Z^{n+1}+\frac{\Z z}{q}.
\end{equation*}
For each integer $t$, the solutions of $a\cdot z\equiv t\pmod q$ form the translate $te_0+L$, where $e_0:=(1,0,\ldots,0)$.
\end{lemma}
\begin{proof}
Since $z_0=1$, the homomorphism $a\mapsto a\cdot z$ from $\Z^{n+1}$ to $\Z/q\Z$ is surjective. Its kernel $L$ has index $q$, so $\det L=q$. Every vector in $\Z^{n+1}+\Z z/q$ pairs integrally with $L$. Moreover, the class of $z/q$ modulo $\Z^{n+1}$ has order exactly $q$. The proposed dual thus has determinant $1/q$, which equals $\det L^*$. The inclusion is an equality. The translate assertion follows by subtracting $te_0$.
\end{proof}

\begin{proposition}\label{prop:inhom}
Let $\xi,\eta\in\Zp$. Let $\lambda,c>0$, and suppose that, at a sufficiently large real $T$, there is no nonzero integral vector $x=(x_0,\ldots,x_n)$ with
\begin{equation*}
\|x\|_\infty\le T,
\qquad
\max_{1\le i\le n}|x_0\xi^i-x_i|_p\le cT^{-1-\lambda}.
\end{equation*}
Then there exists $Q\in\Z[X]$ of degree at most $n$ satisfying
\begin{equation*}
H(Q)\le \frac{n+1}{2}pc^{-1}T^{\lambda},
\qquad
|Q(\xi)-\eta|_p\le cT^{-1-\lambda}.
\end{equation*}
\end{proposition}
\begin{proof}
Choose a power $q=p^h$ such that
\begin{equation*}
c^{-1}T^{1+\lambda}\le q<pc^{-1}T^{1+\lambda},
\end{equation*}
and form $L$ as above. We first claim that
\begin{equation*}
\mu_1(L^*)\ge T/q.
\end{equation*}
Otherwise, a nonzero $y\in L^*$ with $\|y\|_2<T/q$ can be written as $y=b+az/q$, with $b\in\Z^{n+1}$ and $a\in\Z$. Then $x:=qy=qb+az$ is a nonzero integral vector with $\|x\|_\infty<T$. Since $z_0=1$,
\begin{equation*}
x_0\xi^i-x_i=q(b_0\xi^i-b_i)+a(\xi^i-z_i)\in q\Zp
\qquad(1\le i\le n).
\end{equation*}
It follows that its maximum error is at most $q^{-1}\le cT^{-1-\lambda}$, contradicting the hypothesis. This proves the claimed lower bound for $\mu_1(L^*)$.

Choose $t\in\Z$ representing $\eta$ modulo $q\Zp$. By Lemma~\ref{lem:cover}, we clearly have a vector $(a_0,\ldots,a_n)\in te_0+L$ with
\begin{equation*}
\max_i|a_i|\le \frac{(n+1)q}{2T}.
\end{equation*}
For $Q(X):=\sum_{i=0}^n a_iX^i$, we have $Q(\xi)-\eta\in q\Zp$. The claimed estimates now follow from the choice of $q$.
\end{proof}
\medskip

\begin{proof}[Proof of Theorem~\ref{thm:reciprocal}]
At $v=p$, let $\lambda>\hl$. Assume $\xi\in\Zp$. By the definition of $\hl$, the system in Proposition~\ref{prop:inhom} has no nonzero integral solution for arbitrarily large $T$, with $c=1$. We give the construction for any fixed $c>0$ for which this nonexistence condition holds. The height bound in Proposition~\ref{prop:inhom} is uniform in $\eta$, which we choose to enforce the prescribed splitting conditions.\\

Apply Lemma~\ref{lem:CRT} with the set $S\cup\{p\}$. At $p$ choose $r_p\in\Z$ with $\xi-r_p\in p^{2E_p+1}\Zp$ and at every other prime choose $r_\ell=0$. This gives fixed $F$ and
\begin{equation*}
D:=\prod_{\ell\in S\cup\{p\}}\ell^{2E_\ell+1}.
\end{equation*}
Since $G_{r_p}(\xi)\in p^{2E_p+1}\Zp$, $F\equiv G_{r_p}\pmod{p^{2E_p+1}}$ and $D/p^{2E_p+1}$ is a $p$-adic unit, 
\begin{equation*}
\eta:=-F(\xi)/D
\end{equation*}
belongs to $\Zp$. For each such $T$, Proposition~\ref{prop:inhom} provides $Q\in\Z[X]$ of degree at most $n$ with
\begin{equation*}
H(Q)\ll T^\lambda \,\,\, \text{and}\,\,\, |Q(\xi)-\eta|_p\ll T^{-1-\lambda}.
\end{equation*}
Set $P:=F+DQ$. Then
\begin{equation*}
H(P)\ll T^\lambda \,\,\, \text{and}\,\,\,
|P(\xi)|_p=|D|_p|Q(\xi)-\eta|_p\ll T^{-1-\lambda}.
\end{equation*}
Its degree is $n$ and it splits completely at all primes in $S\cup\{p\}$. Let $\alpha\in\Qp$ be the unique root of $P$ satisfying $v_p(\alpha-r_p)\ge E_p+1$, namely the root $\beta_0$ in Proposition~\ref{prop:split}. Since $\xi-r_p\in p^{2E_p+1}\Zp$, we get
\begin{equation*}
0<|\xi-\alpha|_p=\frac{|P(\xi)|_p}{|(n-1)!|_p}\ll T^{-1-\lambda}.
\end{equation*}
The error is nonzero because $\xi$ is not algebraic of degree at most $n$. Since the minimal polynomial of $\alpha$ divides $P$, we have $\deg\alpha\le n$ and $\alpha$ is totally $(S\cup\{p\})$-adic. Now, Lemma~\ref{lem:height} gives
\begin{equation*}
H(\alpha)\ll H(P)\ll T^\lambda.
\end{equation*}
Eliminating $T$ yields
\begin{equation}\label{eq:finalbound}
0<|\xi-\alpha|_p\ll H(\alpha)^{-1-1/\lambda}.
\end{equation}
Indeed, if $H(\alpha)\le A T^\lambda$, then
$T^{-1-\lambda}\le A^{1+1/\lambda}H(\alpha)^{-1-1/\lambda}$. Since $0<|\xi-\alpha|_p\to0$ as $T\to\infty$, there are infinitely many distinct $\alpha$ of unbounded height, as only finitely many algebraic numbers have bounded degree and bounded height.
\end{proof}
\medskip

\begin{proof}[Proof of Theorem~\ref{thm:basic}]
Assume $\xi\in\Zp$. Put $k:=\lceil n/2\rceil$. By Theorem~\ref{thm:BBS}, the system in Proposition~\ref{prop:inhom} has no nonzero integral solution for arbitrarily large $T$ with $\lambda=1/k$ and a fixed $c>0$. The construction in the proof of Theorem~\ref{thm:reciprocal} gives infinitely many distinct algebraic numbers $\alpha\in\Qp$ of degree at most $n$, totally $(S\cup\{p\})$-adic, for which \eqref{eq:finalbound} becomes
\begin{equation*}
0<|\xi-\alpha|_p\ll H(\alpha)^{-1-k}
=H(\alpha)^{-1-\lceil n/2\rceil}.
\end{equation*}
\end{proof}
\medskip

\subsection{The real case}\label{subsec:realproof}

Fix $\xi\in\R$ and $n\ge1$. To obtain totally real approximants, we shall prescribe all the roots of the polynomial while keeping its value at $\xi$ small. We first prove the required inhomogeneous estimate.

\begin{proposition}\label{prop:realinhom}
Let $\lambda>0$, and suppose that for some $T\ge1$ the system
\begin{equation*}
|x_0|\le T,\qquad
\max_{1\le i\le n}|x_i-x_0\xi^i|\le T^{-\lambda}
\end{equation*}
has no nonzero integral solution. Let $D\ge1$ be a fixed integer and $F\in\Z[X]$ a fixed polynomial of degree at most $n$. For every real polynomial $R(X)=\sum_{i=0}^n r_iX^i$ of degree at most $n$, there is $P(X)=\sum_{i=0}^n a_iX^i\in\Z[X]$ with $P\equiv F\pmod D$ and
\begin{equation*}
|P(\xi)-R(\xi)|\le C_0/T,\qquad
|a_i-r_i|\le C_0T^\lambda\quad(1\le i\le n),
\end{equation*}
where $C_0:=D(n+1)/2$.
\end{proposition}
\begin{proof}
Identify polynomials with coefficient vectors in $\R^{n+1}$, and define the invertible real map
\begin{equation*}
A_T(a_0,\ldots,a_n)
:=\bigl(T\textstyle\sum_{i=0}^n a_i\xi^i,
T^{-\lambda}a_1,\ldots,T^{-\lambda}a_n\bigr).
\end{equation*}
For $L_T:=A_T\Z^{n+1}$, the dual lattice is $A_T^{-t}\Z^{n+1}$. Solving $A_T^ty=x$ shows that its vectors have the form
\begin{equation*}
\left(\frac{x_0}{T},\ T^\lambda(x_1-x_0\xi),\ldots,
T^\lambda(x_n-x_0\xi^n)\right),\qquad x\in\Z^{n+1}.
\end{equation*}
If a nonzero vector of $L_T^*$ had Euclidean norm at most $1$, each of its coordinates would have absolute value at most $1$. The corresponding nonzero integral vector $x$ would therefore satisfy
\begin{equation*}
|x_0|\le T,\qquad
\max_{1\le i\le n}|x_i-x_0\xi^i|\le T^{-\lambda},
\end{equation*}
contrary to the hypothesis. Thus $\mu_1(L_T^*)\ge1$, and the covering-radius bound in Lemma~\ref{lem:cover} yields
\begin{equation*}
\rho(L_T)\le\frac{n+1}{2\mu_1(L_T^*)}\le\frac{n+1}{2}.
\end{equation*}
The polynomials congruent to $F$ modulo $D$ have coefficient vectors in $F+D\Z^{n+1}$. By linearity of $A_T$, the image of this set is
\begin{equation*}
A_TF+D L_T.
\end{equation*}
Since $\rho(D L_T)=D\rho(L_T)\le D(n+1)/2=C_0$, there is a vector $b\in\Z^{n+1}$ such that
\begin{equation*}
\|A_T(R-F)-D A_Tb\|_2\le C_0.
\end{equation*}
Let $P$ be the polynomial with coefficient vector $F+Db$. Then $P\in\Z[X]$, $\deg P\le n$ and $P\equiv F\pmod D$. Moreover,
\begin{equation*}
\|A_T(P-R)\|_2\le C_0.
\end{equation*}
Writing $P(X)=\sum_{i=0}^n a_iX^i$, we have
\begin{equation*}
A_T(P-R)=\bigl(T(P(\xi)-R(\xi)),\,
T^{-\lambda}(a_1-r_1),\ldots,T^{-\lambda}(a_n-r_n)\bigr).
\end{equation*}
Each coordinate has absolute value at most $C_0$, so
\begin{equation*}
|P(\xi)-R(\xi)|\le C_0/T,\qquad
|a_i-r_i|\le C_0T^\lambda\quad(1\le i\le n).
\end{equation*}
This proves the proposition.
\end{proof}
\medskip

For real polynomials $U=\sum_{i=0}^n u_iX^i$, we use the same definition of height:
$$
H(U):=\max_{0\le i\le n}|u_i|.
$$
Define the monic polynomial
\begin{equation*}
G(X):=\prod_{j=0}^{n-1}(X-\xi-j).
\end{equation*}
For each $0\le j\le n-1$, let $I_j$ be the open interval
\begin{equation*}
I_j:=(\xi+j-\tfrac14,\xi+j+\tfrac14).
\end{equation*}
These intervals are pairwise disjoint, and each contains exactly one root of $G$. Let $\mathcal E$ denote the set of their $2n$ endpoints. Define the positive constant
\begin{equation*}
\delta:=\frac12\min\left\{1,\ \min_{t\in\mathcal E}
\frac{|G(t)|}{\sum_{i=0}^n|t|^i}\right\}>0.
\end{equation*}
The positivity follows because $G$ has no zero in $\mathcal E$.

\begin{lemma}\label{lem:realstability}
If $U\in\R[X]$ has degree at most $n$ and
$H(U-G)\le\delta$, then $U$ has degree exactly $n$, leading coefficient at least $1/2$, and exactly one simple real root in each $I_j$.
\end{lemma}
\begin{proof}
Since $G$ is monic and $\delta\le1/2$, the leading coefficient of $U$ is at least $1/2$. At any endpoint $t\in\mathcal E$,
\begin{equation*}
|U(t)-G(t)|\le\delta\sum_{i=0}^n|t|^i\le\tfrac12|G(t)|.
\end{equation*}
Thus $U(t)$ and $G(t)$ have the same sign. Since $G$ has one simple root in each $I_j$, its endpoint values have opposite signs. The intermediate value theorem gives a root of $U$ in every $I_j$. There are $n$ disjoint intervals and $\deg U=n$, so these are all its roots and each is simple.
\end{proof}

We now combine Proposition~\ref{prop:realinhom} with Lemma~\ref{lem:realstability} to construct an integral polynomial in a prescribed congruence class, with all its roots real and one close to $\xi$.

\begin{proposition}\label{prop:realroots}
Assume that $\xi$ is not algebraic of degree at most $n$, and fix $D\ge1$ and $F\in\Z[X]$ of degree at most $n$. There are positive constants $C_2,C_3$, depending only on $n,\xi,D$, with the following property. If $\lambda>0$ and $T\ge1$ satisfy the assumption of Proposition~\ref{prop:realinhom}, there is a polynomial $P\in\Z[X]$ of degree exactly $n$ such that
\begin{equation*}
P\equiv F\pmod D,\qquad H(P)\le C_2T^\lambda.
\end{equation*}
It has exactly one simple real root $\alpha_j$ in each interval $I_j$ defined above, and
\begin{equation*}
0<|\xi-\alpha_0|\le C_3T^{-1-\lambda}.
\end{equation*}
\end{proposition}
\begin{proof}
Let $C_0:=D(n+1)/2$, and use $G$ and $\delta$ defined above. Set
\begin{equation*}
C_1:=C_0\sum_{i=0}^n|\xi|^i.
\end{equation*}
The constants $C_0,C_1$ depend only on $n,\xi,D$. We now apply Proposition~\ref{prop:realinhom} with
\begin{equation*}
R_T(X):=(1+C_1/\delta)T^{\lambda}G(X).
\end{equation*}
Since $G(\xi)=0$, the resulting polynomial $P$ satisfies
\begin{equation*}
|P(\xi)|\le C_0/T.
\end{equation*}
Write $P-(1+C_1/\delta)T^{\lambda}G=\sum_{i=0}^n e_iX^i$. For $i\ge1$, Proposition~\ref{prop:realinhom} gives $|e_i|\le C_0T^{\lambda}$. Evaluating at $\xi$, we also get
\begin{equation*}
e_0=P(\xi)-\sum_{i=1}^n e_i\xi^i.
\end{equation*}
As $T\ge1$ and $T^{\lambda}\ge1$, it follows that
\begin{equation*}
|e_0|\le C_0/T+C_0T^{\lambda}\sum_{i=1}^n|\xi|^i\le C_1T^{\lambda}.
\end{equation*}
Thus $|e_i|\le C_1T^{\lambda}$ for every $i$ and hence
\begin{equation*}
H\left(\frac{P}{(1+C_1/\delta)T^{\lambda}}-G\right)
\le\frac{C_1}{1+C_1/\delta}<\delta.
\end{equation*}
By Lemma~\ref{lem:realstability}, $P$ has one simple root in every $I_j$, and its leading coefficient $a_n$ satisfies
\begin{equation*}
a_n\ge (1+C_1/\delta)T^{\lambda}/2.
\end{equation*}
Moreover,
\begin{equation*}
H(P)\le\bigl((1+C_1/\delta)H(G)+C_1\bigr)T^{\lambda}.
\end{equation*}
Thus one may take 
$$
C_2:=(1+C_1/\delta)H(G)+C_1.
$$
For $1\le j\le n-1$, we have $\alpha_j\in I_j=(\xi+j-1/4,\xi+j+1/4)$, so
$|\xi-\alpha_j|\ge j-1/4$. Factoring $P$ over $\R$ and using the bounds $|P(\xi)|\le C_0/T$ and $a_n\ge (1+C_1/\delta)T^{\lambda}/2$, we obtain
\begin{align*}
|\xi-\alpha_0|
&=\frac{|P(\xi)|}{a_n\prod_{j=1}^{n-1}|\xi-\alpha_j|}\\
&\le\frac{2C_0}{(1+C_1/\delta)\prod_{j=1}^{n-1}(j-1/4)}\,T^{-1-\lambda}.
\end{align*}
The factor multiplying $T^{-1-\lambda}$ is a positive constant depending only on $n,\xi,D$, which we denote by $C_3$. When $n=1$, the product is $1$. Finally, $P$ is a nonzero integral polynomial of degree $n$, and $\xi$ is not algebraic of degree at most $n$, so the distance is positive.
\end{proof}
\medskip

\begin{proof}[Proof of Theorem~\ref{thm:reciprocal}]At $v=\infty$, let $\lambda>\widehat\lambda_{n,\infty}(\xi)$. By definition, arbitrarily large $T$ satisfy the assumption of Proposition~\ref{prop:realinhom}. Applying Lemma~\ref{lem:CRT} to $S$, with $r_\ell=0$ at every $\ell\in S$, we obtain fixed $D$ and $F$. For every admissible $T$, Proposition~\ref{prop:realroots} gives a polynomial $P$ with $P\equiv F\pmod D$ and a distinguished real root $\alpha:=\alpha_0$ satisfying
\begin{equation*}
H(P)\le C_2T^\lambda \,\,\, \text{and}\,\,\,
0<|\xi-\alpha|\le C_3T^{-1-\lambda}.
\end{equation*}
The minimal polynomial of $\alpha$ divides $P$. It therefore splits over $\R$ and $\Q_\ell$ for every $\ell\in S$. Thus $\deg\alpha\le n$ and $\alpha$ is totally real and totally $S$-adic. Now, by Lemma~\ref{lem:height},
\begin{equation*}
H(\alpha)\le C_4T^\lambda,
\end{equation*}
where $C_4:=2^n\sqrt{n+1}\,C_2$. Therefore
\begin{equation*}
0<|\xi-\alpha|\le C_3C_4^{1+1/\lambda}
H(\alpha)^{-1-1/\lambda},
\end{equation*}
which proves the bound in Theorem~\ref{thm:reciprocal}. Since $0<|\xi-\alpha|\to0$ as $T\to\infty$, we obtain infinitely many distinct approximants $\alpha$ of unbounded height.
\end{proof}
\medskip

\begin{proof}[Proof of Theorem~\ref{thm:real}]
For $n\ge2$, Theorem~\ref{thm:PRreal} gives
$\widehat\lambda_{n,\infty}(\xi)\le B_n^{-1}$.
Choose $\lambda>B_n^{-1}$ sufficiently close to $B_n^{-1}$ such that
$1/\lambda>B_n-\varepsilon/2$.
Theorem~\ref{thm:reciprocal} provides infinitely many $\alpha$ of degree at most $n$, both totally real and totally $S$-adic, satisfying \eqref{eq:reciprocal} with $v=\infty$.
Their heights tend to infinity along a subsequence, so the constant is absorbed into $H(\alpha)^{\varepsilon/2}$. This proves Theorem~\ref{thm:real}.
\end{proof}

\begin{remark}\label{rem:exactdegree}
One can also require the approximants to have degree exactly $n$. For a fixed prime $q_0\nmid D$, impose the additional congruence $P(X)\equiv X^n+q_0\pmod{q_0^2}$ using the Chinese remainder theorem. Eisenstein's criterion then makes $P$ irreducible, without changing the approximation exponents or prescribed splitting conditions.
\end{remark}

The two constructions give the same error estimate by different local mechanisms: in the real case $|P(\xi)|\ll T^{-1}$ and $|P'(\alpha_0)|\asymp T^\lambda$, whereas in the $p$-adic case, $|P(\xi)|_p\ll T^{-1-\lambda}$ and $|P'(\xi)|_p=|(n-1)!|_p$.
\bigskip

\section{An improvement in even degrees}\label{sec:even}

In even degrees, we obtain the following refinement by adapting the argument of Schleischitz \cite[Theorem~2.3]{Schleischitz}.\\

For $m\ge2$, define the cubic polynomial
\begin{equation*}
f_m(t):=mt^3+(m-1)t^2+m(2m-1)t-(2m-1)
\end{equation*}
and let $\theta_m$ be its unique positive root. As verified below, $1/(m+1)<\theta_m<1/m$.

\begin{theorem}\label{thm:even}
For every $m\ge2$ and every transcendental $\xi\in\Qp$,
\begin{equation*}
\widehat\lambda_{2m,p}(\xi)\le\theta_m.
\end{equation*}
\end{theorem}

Using Theorem~\ref{thm:reciprocal}, we obtain that for every $\varepsilon>0$, there are infinitely many distinct algebraic $\alpha\in\Qp$ of degree at most $2m$, totally $(S\cup\{p\})$-adic, such that
\begin{equation*}
0<|\xi-\alpha|_p\le H(\alpha)^{-1-\theta_m^{-1}+\varepsilon}.
\end{equation*}

For example, $\theta_2^{-1}=2.275279\ldots$ and $\theta_3^{-1}=3.184754\ldots$. In general,
$m<\theta_m^{-1}<m+1/2$ and hence this improves Theorem~\ref{thm:basic}.\\ 

We first recall the exponent inequalities used in the proof. For $j\ge1$, recall that $\lambda_{j,p}(\xi)$ is the supremum of the real numbers $\lambda$ for which, for arbitrarily large $T$, there exists a nonzero vector $(x_0,\ldots,x_j)\in\Z^{j+1}$ satisfying
\begin{equation*}
\max_{0\le i\le j}|x_i|\le T,\qquad
0<\max_{1\le i\le j}|x_0\xi^i-x_i|_p\le T^{-1-\lambda}.
\end{equation*}

For $j\ge1$, let $w_{j,p}(\xi)$ be the supremum of the real numbers $w$ for which, for arbitrarily large $T$, there exists a nonzero polynomial $Q\in\Z[X]$ satisfying
\begin{equation*}
\deg Q\le j,\qquad H(Q)\le T,\qquad
0<|Q(\xi)|_p\le T^{-1-w}.
\end{equation*}
The uniform exponent $\widehat w_{j,p}(\xi)$ is the supremum of the real numbers $w$ for which such a polynomial exists for every sufficiently large $T$. Dirichlet's theorem gives
\begin{equation*}
w_{j,p}(\xi)\ge\widehat w_{j,p}(\xi)\ge j,
\qquad
\lambda_{j,p}(\xi)\ge\widehat\lambda_{j,p}(\xi)\ge1/j.
\end{equation*}

By \cite[Theorem~2.9]{BBS}, \cite[equations~(9) and~(18)]{Schleischitz}, and \cite[Corollary~10.3]{PR}, we have the following bounds:
\begin{align}
\widehat\lambda_{a+b-1,p}(\xi)
&\le\max\{w_{a,p}(\xi)^{-1},\widehat w_{b,p}(\xi)^{-1}\},
&&a,b\ge1,\label{eq:cross}\\
\lambda_{j,p}(\xi)&\le\frac{w_{j,p}(\xi)-j+1}{j},
&&j\ge1,\label{eq:khintchine}\\
\lambda_{N,p}(\xi)&\ge
\frac{\widehat\lambda_{N,p}(\xi)^2+(N-2)\widehat\lambda_{N,p}(\xi)}
{(N-1)(1-\widehat\lambda_{N,p}(\xi))},
&&N\ge2,\label{eq:SS}
\end{align}
where \eqref{eq:SS} is used only when $\widehat\lambda_{N,p}(\xi)<1$, and $1/\infty=0$.
\medskip

\begin{proof}[Proof of Theorem~\ref{thm:even}]
By Theorem~\ref{thm:BBS}, $0<\widehat\lambda_{2m,p}(\xi)\le1/m<1$. Applying \eqref{eq:cross} with $a=m$, $b=m+1$, and using $\widehat w_{m+1,p}(\xi)\ge m+1$, we obtain
\begin{equation}\label{eq:firstbalance}
\widehat\lambda_{2m,p}(\xi)\le\max\{1/w_{m,p}(\xi),1/(m+1)\}.
\end{equation}
Restricting to the first $m+1$ coordinates gives $\lambda_{2m,p}(\xi)\le\lambda_{m,p}(\xi)$. Indeed, a vector with $x_0=\cdots=x_m=0$ and a remaining nonzero coordinate would satisfy
\begin{equation*}
|x_i|_p\ge |x_i|^{-1}\ge T^{-1},
\end{equation*}
which is incompatible with a positive simultaneous exponent for large $T$. Combining coordinate restriction with \eqref{eq:khintchine} and \eqref{eq:SS}, we get
\begin{equation}\label{eq:secondbalance}
\frac{\widehat\lambda_{2m,p}(\xi)^2+(2m-2)\widehat\lambda_{2m,p}(\xi)}{(2m-1)(1-\widehat\lambda_{2m,p}(\xi))}
\le\lambda_{2m,p}(\xi)\le\lambda_{m,p}(\xi)\le\frac{w_{m,p}(\xi)-m+1}{m}.
\end{equation}

The polynomial $f_m$ is strictly increasing for $t\ge0$, because
\begin{equation*}
f_m'(t)=3mt^2+2(m-1)t+m(2m-1)>0.
\end{equation*}
Moreover,
\begin{equation*}
f_m\left(\frac1{m+1}\right)
=-\frac{m(2m^2+2m-1)}{(m+1)^3}<0,
\qquad
f_m\left(\frac1m\right)=\frac1m>0.
\end{equation*}
Thus its unique positive root $\theta_m$ lies in $(1/(m+1),1/m)$.
Moreover,
\begin{equation*}
f_m\left(\frac{2}{2m+1}\right)
=-\frac{(2m-1)(4m^2-3)}{(2m+1)^3}<0
\end{equation*}
and hence $\theta_m^{-1}<m+1/2$.\\

If $\widehat\lambda_{2m,p}(\xi)\le1/(m+1)$, then $\widehat\lambda_{2m,p}(\xi)<\theta_m$. Otherwise, \eqref{eq:firstbalance} gives $w_{m,p}(\xi)\le1/\widehat\lambda_{2m,p}(\xi)$, so $w_{m,p}(\xi)$ is finite. Equation~\eqref{eq:secondbalance} yields
\begin{equation*}
\frac{\widehat\lambda_{2m,p}(\xi)^2+(2m-2)\widehat\lambda_{2m,p}(\xi)}{(2m-1)(1-\widehat\lambda_{2m,p}(\xi))}
\le\frac{\widehat\lambda_{2m,p}(\xi)^{-1}-m+1}{m}.
\end{equation*}
Multiplying by the positive number $m(2m-1)\widehat\lambda_{2m,p}(\xi)(1-\widehat\lambda_{2m,p}(\xi))$ and collecting terms gives $f_m(\widehat\lambda_{2m,p}(\xi))\le0$. Hence $\widehat\lambda_{2m,p}(\xi)\le\theta_m$, as required.
\end{proof}

\bigskip

\section{Concluding remarks}\label{sec:remarks}

A natural question is whether our main results extend from $\Q$ to a number field $K$. In this setting, we seek to approximate $\xi\in K_v$ by algebraic numbers $\alpha$ of bounded degree over $K$, whose minimal polynomials split completely over $K_v$ and at a prescribed finite set of other places. We measure the approximation error using the absolute value at $v$ and compare it with a suitable height of $\alpha$. The adelic transference methods of \cite{Roy,RW} may be useful here.\\

Our results establish Wirsing-type approximation where the approximants are restricted to being totally real or totally $p$-adic. It may, however, be premature to expect the full Wirsing's conjecture to be attainable under these restrictions, namely, that for $\xi\in\Q_v$ transcendental and any $\varepsilon > 0$, there exist infinitely many algebraic numbers $\alpha$ of degree at most $n$, totally $v$-adic such that
$$
    0<|\xi - \alpha|_v \leq H(\alpha)^{-n-1-\varepsilon}.
$$
A desirable intermediate step would be to establish analogues of Po\"{e}ls's bound with approximants restricted to these extensions. At present, we do not know how to obtain such a result.\\

Approximation by algebraic numbers belonging to a fixed infinite extension deserves systematic study. The quality of approximation attainable within such an extension may provide information about its arithmetic structure. For instance, we demonstrated that the quality of approximation by $\Q^{(2)}$ is far inferior to that of $\Q^{tr}$ or $\Q^{tp}$. This would suggest that there are ``fewer'' elements in $\Q^{(2)}$, which is consistent with the fact that $\Q^{(2)}$ satisfies the Northcott property with respect to the Weil height, where as $\Q^{tr}$ or $\Q^{tp}$ do not. More generally, one may ask how does this quality of approximation depends on an infinite extension. We relegate this to future investigation.

\bigskip

\section*{Declaration of AI use}
The authors used GPT-6-Astra as an exploratory tool and to identify and correct mathematical errors in certain arguments. In particular, the tool suggested using Schleischitz's result to improve the bound in the even-degree case.

\renewcommand{\bibliofont}{\small}

\end{document}